\documentclass[12pt]{amsart}
\usepackage{amsmath}
\usepackage{amsxtra}
\usepackage{amscd}
\usepackage{amsthm}
\usepackage{amsfonts}
\usepackage{amssymb}

\newtheorem{theorem}{Theorem}[section]

\theoremstyle{definition}

\theoremstyle{plain}
\newtheorem{lemma}[theorem]{Lemma}

\theoremstyle{remark}

\newenvironment{example}
{\noindent{\bf Example\/}:}{\par} \numberwithin{equation}{section}

\newcommand{\U}{{\mathbb U}}
\newcommand{\W}{{\mathbb W}}
\newcommand{\X}{{\mathbb X}}

\newcommand{\R}{{\mathbb R}}
\newcommand{\C}{{\mathbb C}}

\newcommand{\Pee}{{\mathbb P}}

\newcommand{\la}{\lambda}

\newcommand{\HH}{\mathcal{H}}

\newcommand{\al}{\alpha}
\newcommand{\bean}{\begin{eqnarray}}
\newcommand{\eean}{\end{eqnarray}}
\newcommand{\be}{\begin{displaymath}}
\newcommand{\ee}{\end{displaymath}}
\newcommand{\bea}{\begin{eqnarray*}}
\newcommand{\eea}{\end{eqnarray*}}

\newcommand{\Gr}{{\rm Gr}}

\begin{document}

\title[The period matrices of an arrangement and its dual]
{Determinants of the hypergeometric
period matrices of an arrangement and
its dual}

\author[D. Mukherjee and A. Varchenko]
{D. Mukherjee and A. Varchenko ${}^{1}$}

\thanks{${}^1$ Supported in part by NSF grant DMS-0244579}

\begin{abstract}
We fix three natural numbers $k, n, N$, such that
$n+k+1=N$, and introduce the
notion of two dual arrangements of hyperplanes. One of the
arrangements is an arrangement of $N$ hyperplanes in a
$k$-dimensional affine space, the other is an arrangement of $N$
hyperplanes in an $n$-dimensional affine space. We assign weights
$\al_1, \dots, \al_N$ to the hyperplanes of the arrangements and for
each of the arrangements consider the associated period matrices.
The first is a matrix of $k$-dimensional hypergeometric integrals
and the second is a matrix of $n$-dimensional hypergeometric
integrals. The size of each matrix is equal to the number of bounded
domains of the corresponding arrangement.  We show that the dual
arrangements have the same number of bounded domains and the product
of the determinants of the period matrices is equal to an
alternating product of certain values of Euler's gamma function
multiplied by a product of exponentials of the weights.

\end{abstract}

\maketitle \centerline{\it Department of Mathematics, University of
  North Carolina at Chapel Hill,} \centerline{\it Chapel Hill, NC
  27599-3250, USA} \medskip

\section{Introduction}

Let $\mathcal{V}$ be a real affine space of dimension $k$. Let
$\mathcal{F} = \{\,f^j\ :\ j \in \mathcal{J} \,\}$ be a finite set of degree
one polynomials defined on $\mathcal{V}$. For $j \in \mathcal{J}$, let
$\mathcal{H}^j$ be the hyperplane in $\mathcal{V}$ given by the
zero-set of $f^j$. Consider the affine hyperplane arrangement
$\mathcal{A} = \{\,\mathcal{H}^j\ :\ j \in \mathcal{J} \,\}$.  Assume that a
positive number $\alpha_j$ is assigned to each
hyperplane $\mathcal{H}^j$.

A bounded connected component of $\mathcal{V} \, \setminus \, \cup_{j \in \mathcal{J}}
\mathcal{H}^j $ is called {\it a bounded domain} of
$\mathcal{A}$.
Let $\text{\tt Ch}(\mathcal{A})$ be the set of bounded domains of
$\mathcal{A}$ and $\beta$ the number of bounded domains.

A { logarithmic differential $k$-form} associated to
$\mathcal{F}$ is a $k$-form of the type
$$ \phi =
\sum a_{j_1 \, j_2 \, \ldots \, j_k } \frac{df^{j_1}}{f^{j_1}}
\wedge \frac{df^{j_2}}{f^{j_2}} \wedge \cdots \wedge
\frac{df^{j_k}}{f^{j_k}}
$$
with $ a_{j_1 \, j_2 \, \ldots \,j_k } \in \R$.
The form is regular on $\mathcal{V} \setminus \cup_{j \in
\mathcal{J}} \mathcal{H}^j$.

If $\mathcal{F}$ is ordered, then, using constructions from
 \cite{DT} and \cite{V1}, one obtains
\begin{enumerate}
\item[$\bullet$] an ordered set $\Phi = \{ \phi^1, \phi^2 , \ldots , \phi^{\beta}\}$
of logarithmic differential $k$-forms,
\item[$\bullet$] an order on $\text{\tt Ch}(\mathcal{A}) =
\{\Delta_1 , \Delta_2 , \ldots , \Delta_{\beta} \}$,
\item[$\bullet$] an orientation of each $\Delta \in \text{\tt Ch}(\mathcal{A})$.
\end{enumerate}

Consider the multi-valued function
$$
U^{\alpha} = \prod_{j \in \mathcal{J}} (f^j)^{\alpha_j} .
$$
Fix a uni-valued branch of each function
$(f^j)^{\alpha_j}$ on each bounded domain of
the arrangement. This determines the uni-valued branch of
the function $U^{\alpha}$ on each
bounded domain. The $\beta \times \beta$-matrix
$$
\text{PM}\, ( \mathcal{A} , \alpha )  = \left( \int_{\Delta_s}
U^{\alpha} \phi^t \right) , \qquad s,t= 1 , 2 , \ldots , \beta ,
$$
is called  {\it the period matrix of the weighted ordered arrangement}.

\medskip

\noindent
{\bf Example 1.}
Let $a_1 < a_2 < a_3$ be real numbers. Consider the set of
polynomials
$$
\mathcal{F} = \{\,f^j\,=\,x-a_j \ :\ j=1,2,3 \,\}
$$
defined on
$\R$. Then the arrangement $\mathcal{A}= \{ \mathcal{H}^1,
\mathcal{H}^2, \mathcal{H}^3 \}$ is the set of three points $a_1 ,
a_2 , a_3$ in $\R$. Let $\alpha_j$ be the weight of $\mathcal{H}^j$, $j =
1,2,3$.

There are two bounded domains:\ $ \Delta_1 = (a_1 , a_2),    \Delta_2
= (a_2 , a_3)$. The set of $1$-forms is
$$
\Phi\ =\ \{\ \phi^1 = \alpha_2 \frac{dx}{x-a_2}\ ,\;\; \phi^2 = \alpha_3
\frac{dx}{x-a_3}\ \}\ .
$$
Consider the function
$$
U^{\alpha}\ = \ (x-a_1)^{\alpha_1}
(x-a_2)^{\alpha_2} (x-a_3)^{\alpha_3} \ .
$$
Fix a uni-valued branch of each function $(x-a_j)^{\alpha_j}$
on each interval $\Delta_s$. The period matrix is
$$
\left(
\begin{array}{cc}
\int_{a_1}^{a_2} \alpha_2 \prod_{j=1}^3 (x-a_j)^{\alpha_j}
\frac{dx}{x-a_2} &
\int_{a_1}^{a_2} \alpha_3 \prod_{j=1}^3 (x-a_j)^{\alpha_j} \frac{dx}{x-a_3} \\

\noalign{\bigskip}

\noalign{\medskip}

\int_{a_2}^{a_3} \alpha_2 \prod_{j=1}^3 (x-a_j)^{\alpha_j}
\frac{dx}{x-a_2} & \int_{a_2}^{a_3} \alpha_3 \prod_{j=1}^3
(x-a_j)^{\alpha_j} \frac{dx}{x-a_3}
\end{array}
\right) .
$$

\medskip

In \cite{V1} and \cite{DT}, the determinant of the period matrix was
computed in terms of {\it critical values}\, $c \,((f^j)^{\alpha_j} ,
\Delta )$ of the chosen branches of the functions $(f^j)^{\alpha_j}$
on the bounded domains and a certain function, called {\it the beta
function of the weighted arrangement}, see Sections
\ref{criticalvalues} and \ref{betafunction}. The beta function is an
alternating product of values of Euler's gamma function whose arguments
are appropriate linear combinations of the $\alpha_j$'s. It is proved
in \cite{V1} and \cite{DT}, that the determinant of the period matrix
is given by the formula:
$$
\det \left( \int_{\Delta_s} U^{\alpha} \phi^t \right)\ =\
{\rm B} ( \mathcal{A} , \alpha ) \,\cdot
 \prod_{\stackrel{\Delta \in \;
\text{\tt Ch}(\mathcal{A})}{j \in \mathcal{J}}} c
\,((f^j)^{\alpha_j} , \Delta )\ ,
$$
where ${\rm B} ( \mathcal{A} , \alpha )$
is the beta function of the weighted arrangement.

\medskip

\noindent
{\bf Example 2.}
The beta function of the arrangement  in Example 1 is
$$
{\rm B} (\mathcal{A} , \alpha) = \frac{\Gamma (\alpha_1 +1)
\Gamma(\alpha_2 +1) \Gamma (\alpha_3 +1)}{\Gamma (\alpha_1 +
\alpha_2 + \alpha_3 +1)}\ .
$$
The product of the critical values is
$$
(f^1)^{\alpha_1}_{\Delta_1} (a_2) \cdot (f^1)^{\alpha_1}_{\Delta_2}
(a_3) \cdot (f^2)^{\alpha_2}_{\Delta_1} (a_1) \cdot
(f^2)^{\alpha_2}_{\Delta_2} (a_3) \cdot (f^3)^{\alpha_3}_{\Delta_1}
(a_1) \cdot (f^3)^{\alpha_3}_{\Delta_2} (a_2)\ ,
$$
where $(f^j)^{\alpha_j}_{\Delta_s} $ is the chosen branch of
$(f^j)^{\alpha_j}$ on $\Delta_s$, and $(f^j)^{\alpha_j}_{\Delta_s}(a_i)$ is the value
of that branch at $a_i$.

\medskip

In this paper, we introduce the notion of dual arrangements.
We fix natural numbers $k, n, N$ such that
\bea
k+n+1=N\ ,
\qquad
3\leq N\ ,
\qquad
1 \leq k,n \leq N-2\ ,
\eea
and consider the vector space $ \R^{N+1}$ and its dual space.
Let $\{e_1 ,  \ldots ,$ $ e_{N+1} \}$ be the standard basis of $\R^{N+1}$
and $\{ e^1 , \ldots , e^{N+1} \}$ the dual basis of the dual space.
We denote $\R^{N+1}$ by $\X$ and the dual space by $\X^{\prime}$.
Set $J = \{ 1 , \ldots , N, N+1 \}$ and $\mathcal J = \{ 1 , \ldots , N \}$.

Let ${\mathbb W} \subset \X$ be a vector subspace of dimension
$k+1$. Let $\W'\subset \X'$ be the annihilator of $\W$.  The
subspace $\W'$  is of dimension $n+1$.
We assume that for any $a,b \in J$, $a\neq b$, the functions $e^a
\vert_{\W}$ and $e^b \vert_{\W}$ are not proportional, and the
functions $e_a \vert_{\W'}$ and $e_b \vert_{\W'}$ are not
proportional.

The pair $\tau = (\X , \W )$ with this property will be called {\it
an admissible pair in $\X$}.  Similarly, the pair
$\tau' = (\X' , \W' )$ with this property will be called {\it an
admissible pair in $\X'$}. The pairs $\tau $ and
$\tau'$ will be called {\it dual}.

Let $V = \Pee (\W)$ be the projective space associated with $\W$
and $\mathcal{V} \subset V$ the affine space defined by
the condition $e^{N+1}\vert_\W \neq 0$. The spaces $V$ and
$\mathcal V$ are of dimension $k$.
The functions $e^j/e^{N+1}$, $ j \in \mathcal J$, define on $\mathcal V$ a set of degree
one polynomials $f^j$ and an arrangement of hyperplanes denoted by
$\mathcal A[\tau]$.

Similarly, let $V' = \Pee (\W')$ be the projective space associated with $\W'$
and $\mathcal{V}' \subset V'$ the affine space defined by
the condition $e_{N+1}\vert_{\W'} \neq 0$. These are spaces of dimension $n$.
The functions $e_j/e_{N+1}$, $ j \in \mathcal J$, define on $\mathcal V'$ a set of degree
one polynomials $f_j$ and an arrangement of hyperplanes denoted by
$\mathcal A[\tau']$.

The arrangements $\mathcal A[\tau]$ and $\mathcal A[\tau']$ will be called
{\it dual}. These are arrangements of $N$ hyperplanes in affine
spaces of dimensions $k$ and $n$, respectively.  We prove that the
number of bounded domains of $\mathcal A[\tau]$ is equal to the number of
bounded domains of $\mathcal A[\tau']$ and study other combinatorial similarities
between the dual arrangements.

Fix positive numbers $\{ \alpha_j :  j \in \mathcal J \}$,
then arrangements $\mathcal A[\tau]$ and $\mathcal A[\tau']$ become
weighted arrangements.  Let $\text{PM}\, ( \mathcal{A}[\tau] , \alpha
)$ and $\text{PM}\, ( \mathcal{A}[\tau'] , \alpha )$ be their period
matrices. The period matrices depend on the choice of uni-valued
branches of the corresponding functions on the corresponding bounded
domains. In this paper, we give a construction of the choice of
branches so that the determinants of the period matrices become
related. For our choice of branches we prove that
\bea
\det\, \text{PM}\, ( \mathcal{A}[\tau] , \alpha )\ \cdot\
\det \, \text{PM}\, ( \mathcal{A}[\tau'] , \alpha ) \ = \
\left[ \frac{\prod_{j \in \mathcal{J}}\, e^{\pi i \alpha_j}\,
\Gamma (\alpha_j + 1)}{\Gamma ( \sum_{ j \in \mathcal{J}} \alpha_j + 1 )}
\right]^{\beta} \,  ,
\eea
see Theorem \ref{main theorem}. This formula relates the determinants
of matrices of $k$-dimensional and $n$-dimensional hypergeometric
integrals and shows that the product of determinants of period
matrices of dual arrangements is a combinatorial quantity, not
depending on the size of the bounded domains or angles between
hyperplanes.  Theorem \ref{main theorem} is the main result of the
paper.

\medskip

\noindent
{\bf Example 3.}
For the arrangement in Example 1, the dual
arrangement is also an arrangement of three points on the real
line. The points $\mathcal{H}_1 , \mathcal{H}_2 , \mathcal{H}_3$ of
the dual arrangement are given by the zero-sets of the polynomials
$$
f_1 = \frac{a_2 -a_3}{a_1 - a_2} \left( x- \frac{1}{a_3 - a_2}
\right), \quad f_2 = \frac{a_3 -a_1}{a_1 - a_2} \left( x-
\frac{1}{a_3 - a_1} \right), \quad  f_3 = x \ ,
$$
respectively. According to Theorem \ref{main theorem}, the product of determinants of
the period matrices of the arrangement of Example 1 and its dual
 is equal to
$$
e^{2 \pi i\, (\alpha_1 + \alpha_2 + \alpha_3) }  \left[ \frac{\Gamma
(\alpha_1 + 1)\, \Gamma (\alpha_2 + 1)\, \Gamma (\alpha_3 + 1)}
{\Gamma(\alpha_1 +\alpha_2 +\alpha_3 + 1)}\right]^{2} \ .
$$

The paper has the following structure. In Section \ref{project}
we discuss combinatorics of an arrangement of hyperplanes.
In Section \ref{dualarrangement} we introduce the notion of dual arrangements
and compare the combinatorics of dual arrangements. Section
\ref{Sec DET} is about period matrices. The section contains the statement
of the main result of the paper, Theorem \ref{main theorem}. In Section
\ref{proofs} we prove Theorem \ref{main theorem}. In Appendix A, we introduce
the notion of weak duality and show that natural constructions with dual arrangements lead
to weakly dual arrangements. In Appendix B, we formulate a statement 
 which helps to
determine if two given arrangements are dual.

\bigskip
The authors thank E.M. Rains who referred the second author to the
paper by A.L. Dixon \cite{D}, published in 1905, in which certain
hypergeometric integrals of different dimensions were equated; see
an elliptic version of Dixon's identity in \cite{R}. E.M. Rains
suggested that there might be similar identities for determinants of
periods of suitable arrangements of different dimensions. It is not yet
clear to us if Dixon's result is related to our Theorem
\ref{main theorem}.

The authors thank T.A. Brylawski for teaching the basics of the
matroid theory.

\section{Arrangements and Matroids}\label{project}

All vector and affine spaces in this paper are over the field of
real numbers.  For a vector space $U$, \ $\Pee (U)$ denotes the
projective space of one-dimensional vector subspaces of $U$.

\subsection{Arrangement and edges} \label{arrangement and edges}

Let $\W$ be a vector space and $\Sigma = \{u^j: j \in J\}$ a finite
collection of nonzero vectors in the dual space of $\W$. For $j \in
J$, denote by $E^j\subset \W$ the hyperplane $\{z : u^j (z) =0 \}$
and by  $H^j = \Pee (E^j) \subset \Pee (\W)$ its projectivization.
Then ${A} = \{ {H}^j :  j \in J\}$
 is an  arrangement of hyperplanes in $\Pee (\W)$.

 A non-empty intersection of some of the hyperplanes of the
arrangement is called an {\it edge}. A {\it vertex} is a
zero-dimensional edge.

To an edge $L$, we associate two arrangements:
\begin{enumerate}
\item[$\bullet$]
 ${A}^L = \{ {H}^j : L \subset {H}^j \}$,
{ \it the localization of ${A}$ at $L$},

\item[$\bullet$]
 ${A}_L = \{{H}^j \cap L : L \not
\subset {H}^j \}$, { \it the induced arrangement on $L$}.
\end{enumerate}

An arrangement ${A} = \{ {H}^j :  j \in J\}$
is called {\it central} if the intersection $L=
\cap_{j \in J} {H}^j$ is not empty.

For a central arrangement, consider the projective space $\Pee_L$
whose points are the $(\dim L + 1)$-dimensional projective subspaces
of $\Pee(\W)$ containing $L$. Any hyperplane in $\Pee(\W)$,
containing $L$, determines uniquely a hyperplane in $\Pee_L$. Thus,
the central arrangement determines an arrangement of hyperplanes in
$\Pee_L$ called {\it the projectivization of the central
arrangement.}

The projectivization of the central arrangement ${A}^L$ is
 called {\it the projective localization of} $L$
and denoted by $P({A}^L)$.

\medskip

Let $L$ be an edge. Let $a,b\in J$, ${a}\neq {b}$.  We say that $L$
is parallel to $ H^{a}$ in the affine space $ V \setminus H^{b}$, if
$H^a$ does not coincide with $H^b$ and $L$ does not intersect $
H^{a}$ in the affine space.  In that case we also say that the
triple $(L,H^{a},H^{b})$ is {\it a parallelism in the arrangement
$A$.}

\subsection{Matroids.} \label{arrandmat}
On matroid theory see \cite{O} and \cite{B}.

Let $J$ be a finite set and $\mathcal{I}$ a collection of subsets of
$J$. The pair $M = ( J , \mathcal{I}\, )$ is called a {\it matroid}
if the following properties hold.
\begin{enumerate}
\item[$\mathcal{I}$1.]  $\emptyset \in \mathcal{I}$.

\item[$\mathcal{I}$2.] If $X$ is in $\mathcal{I}$ and $Y \subset X$ , then $Y$
is also in $\mathcal{I}$.

\item[$\mathcal{I}$3.] If $X$ and $Y$ are in $\mathcal{I}$ and $|X| > |Y|$,
then there is an element $ x \in X \setminus Y$ such that $ Y \cup
\{x\} $ is in $\mathcal{I}$.
\end{enumerate}

The set $J$ is called {\it the ground set} of the matroid and
elements of $\mathcal{I}$ are called  {\it the independent sets} of
the matroid.
\medskip

\begin{example}
\label{vectormatroid} Let $\Sigma = \{ u^j : j \in {J}\}$ be a
finite collection of vectors in a vector space. Define the
collection $\mathcal{I}$ of subsets of $J$: a subset $S \subset J$
belongs to $\mathcal{I}$ if and only if the vectors $\{ u^j : j \in
S \}$ are linearly independent. This defines {\it the matroid of the
collection of vectors}.
\end{example}

\medskip

\begin{example}
\label{projectivematroid} Let $\W$ be a vector space. Let $\Sigma =
\{ u^j : j \in J\}$ be a finite collection of nonzero vectors in the
dual space of $\W$.  The collection defines the  arrangement ${A}$
of hyperplanes in $\Pee (\W)$. The matroid of $\Sigma$ is called
{\it the matroid of the arrangement}.

The arrangement  defines vectors of $\Sigma$ up to multiplication by
nonzero numbers. This multiplication does not change the matroid of
$\Sigma$. Hence, the matroid of the arrangement does not depend on
the choice of the collection of vectors.
\end{example}

\medskip

Let $M = (J, \mathcal{I}\,)$ be a matroid. A maximal (with respect
to inclusion) element of $\mathcal{I}$ is called a {\it basis}. The
axiom ($\mathcal{I}$3) implies that all bases have the same
cardinality. More generally, for any subset $X$ of $J$, the maximal
independent subsets of $X$ all have the same cardinality. Define
\begin{enumerate}
\item[$\bullet$]
$\text{ rank}_M  X $ to be the cardinality of the largest
independent subset of $X$,

\item[$\bullet$] $\text{ corank}_M X =  \text{rank}_M  J  - \text{rank}_M  X $,

\item [$\bullet$]
$\text{ nullity}_M  X  = |X| - \text{rank}_M  X $,

\item [$\bullet$]
${}$ $\hbox{rank } M=\text{ rank}_M  J$.
\end{enumerate}

\subsubsection{Tutte polynomial} The {\it Tutte polynomial} of a
matroid $M = (J , \mathcal{I} \,)$ is the polynomial in $x$ and $y$,
given by the formula
$$
{\rm T} (M; x, y) = \sum_{X \subset J} (x-1)^{\text{corank}_M X} \,
(y-1)^{\text{nullity}_M X} .
$$

\begin{theorem} [\cite{B}] \label{b10}
Let $M$ be a matroid on the ground set $J$. Let ${\rm T} (M; x, y) =
\sum_{i,j} b^{ij}_{M}\, x^i y^j$.
\begin{enumerate}
\item[$\bullet$]
 If $|J| \geq 2$, then $b^{1
0}_{M} = b^{0 1}_{M}$. \\
\item[$\bullet$]
 If $|J| \geq 1$, then $b^{00}_{M} = 0$.
\end{enumerate}
\end{theorem}

\subsubsection{Contraction and deletion} Let $M=(J , \mathcal{I}\,)$ be
a matroid. For a subset $X \subset J$, denote by $\hat{X} = J
\setminus X$ its complement.

For a non-empty subset $X\subset J$, $|X|<|J|$, define the matroid
$M / X = (\hat{X} , \mathcal{I}_{M/X})$ called the {\it contraction
of $X$}. A subset $I \subset \hat{X}$ is in $\mathcal{I}_{M/X}$ if
and only if for some maximal independent subset $Y$ of $X$ in $M$,
the set $I \cup Y$ is independent in $M$.

For a non-empty subset $X \subset J$, $|X|<|J|$, define the matroid
$M - X = (\hat X, \mathcal{I}_{M - X})$ called the {\it deletion of
$X$}. A subset $I \subset \hat X$ is in $\mathcal{I}_{M - X}$ if and
only if $I$ is independent in $M$.

An element $j \in J$ is called a {\it loop} if it is not contained
in any basis of $M$. Dually, an element $j$ is called an {\it
isthmus} if it is contained in every basis.

\begin{theorem} [\cite{B}] \label{tuttepolynomial}
If $j$ is neither a loop nor an isthmus, then
$${\rm T} (M; x, y) = {\rm T} (M - \{j\} ; x, y) + {\rm T} (M /\{j\}  ; x, y).$$
If $j$ is an isthmus, then
$$
{\rm T} (M; x, y) = x \, {\rm T} (M / \{j\}  ; x, y).
$$
\end{theorem}

A non-empty subset $X \subset J$ is called a {\it flat}  if for
every $y \in J \setminus X$, $$\hbox{rank}_M X \cup \{ y \}
> \hbox{rank}_M X .$$

For a flat $X$, define its {\it discrete length, width, volume} as
the numbers \bea {\rm l}_M X \, = \,b^{10}_{M/X}\ ,\, \qquad {\rm
w}_M X \,=\, b^{10}_{M - \hat{X}}\ , \qquad {\rm vol}_M X \,=\,{\rm
l}_M X \cdot {\rm w}_M  X \ , \eea respectively. We say that a  flat
is {\it spacious} if it has a nonzero discrete volume.

Let $X$ be a flat. Let $a, b \in J$, $a \neq b$.  The triple $(X, a,
b)$ is called a {\it parallelism} in $M$ if
 $a, b \not \in X$, ${\rm rank}_M \{ a,b\} =2$, and
 $\text{rank}_M X \cup \{ a  , b \} = \text{rank}_M X + 1$.

Let $(X , a , b)$ be a parallelism. Denote by $\hat X(a,b) = J
\setminus (X\cup \{a,b\})$ the complement of $X \cup \{a,b \}$ in
$J$.

For a parallelism $(X , a , b)$, define its {\it discrete
 width},
{\it volume} as the numbers
$$
{\rm w}_M (X , a , b) = b^{10}_{M - \hat X(a,b)}\ , \qquad {\rm
vol}_M (X , a , b) = {\rm l}_M X \cdot {\rm w}_M (X , a , b)\ ,
$$
respectively.

\subsection{Matroid of an arrangement}

Let $M$ be the matroid of an arrangement \linebreak ${A} = \{ H^j :
j\in  J\}$ of hyperplanes  in a projective space ${V}$.  The flats
in $M$ are in one-to-one correspondence with edges of ${A}$. If $L$
is an edge, then $X = \{j: L \subset {H}^j \}$ is a flat.

Let $L$ be an edge and $X$ the corresponding flat.  Let ${A}^{L}$
and ${A}_{L}$ be the localization and  induced arrangements,
respectively.  Let $M[{A}^{L}]$ and $M[{A}_{L}]$ be the matroids
associated to the arrangements ${A}^{L}$ and ${A}_{L}$,
respectively. Then
$$
 M[{A}^{L}]  \ = \ M  - \hat{X}\ ,
\qquad
 M[{A}_{L}] \ = \ M  / X\ .
$$

Let $P ( {A}^L )$ be the projective localization of $L$. Then the
matroid of ${A}^L$ is also the matroid of $P (A^L)$.

\begin{lemma} \label{fullrank}
If $\hbox{\rm rank } M = \dim {V} + 1$, then $ \hbox{\rm rank } M /
X =  \dim L + 1 . $ \hfill $\square$
\end{lemma}

\subsection{Edges and parallelisms}

Define  {\it the discrete length, width, and volume} of an edge $L$
as the discrete length, width, and volume, respectively, of the flat
$X$, \bea {\rm l}_{{ A}} L \, = \,b^{10}_{M[{A}_{L}]}\ ,\, \qquad
{\rm w}_{{A}}L \,=\, b^{10}_{ M[{A}^{L}]}\ , \qquad {\rm vol}_{{ A}}
L \,=\, {\rm l}_{{ A}} L \cdot {\rm w}_{A} L \, , \eea cf. \cite
{V1}. An edge will be called {\it spacious} if it has a nonzero
discrete volume.

Let $L$ be an edge and $a, b \in J$. The edge $L$ is parallel to $
H^{a}$ in $ V \setminus  H^{b}$, if and only if the triple $(X,a,b)$
is a parallelism in the matroid $M$.

Define  {\it the discrete  width and volume of a parallelism
$(L,H^{a},H^{b})$} in $A$ as the discrete width and volume of the
parallelism $(X , a , b)$ in $M$. That is,
$$
{\rm w}_{A} (L , H^{a} , H^{b}) = {\rm w}_M(X,a,b)\, , \qquad {\rm
vol}_{A} (L , H^{a} , H^{b}) = {\rm l}_A L \cdot {\rm w}_{A} (L ,
H^{a} , H^{b}) .
$$

\subsection{Bounded domains}
Let $A = \{ H^j : j \in J \}$ be an arrangement of hyperplanes in a
projective space $V$. The connected components of the topological
space \linebreak ${V}\, \setminus\, \cup_{j \in J} {H}^j$ are called
{\it domains}. For $j\in J$, a domain is called {\it bounded with
respect to the hyperplane} ${H}^j$ if the closure of the domain does
not intersect the hyperplane.

\begin{theorem} [\cite{ZZ}] ${}$ \label{infinity}
 Assume that $\hbox{\rm rank } M = \dim \, V\, +\, 1$. Then
for $j \in J$, the number of domains of ${A}$ bounded with respect
to ${H}^{j}$ is equal to $b^{10}_{M}$. In particular, the number of
bounded domains does not depend on the choice of $j$.
\end{theorem}

If   $\hbox{\rm rank } M = \dim \, V\, +\, 1$, then
 the discrete length and width of an edge $L$
are the numbers of bounded domains in arrangements ${{ A}}_L$ and
$P(A^L)$, respectively. See Lemma \ref{fullrank}, Theorem
\ref{infinity}.

\subsection{ Geometric interpretation of the
discrete volume of a parallelism in the arrangement $A$} Assume that
$H^{a}\neq H^{b}$.  In the affine space $ \mathcal{V} = V \setminus
H^{b}$, consider the arrangement of hyperplanes $\mathcal A=
\{\,\mathcal{H}^{j}\ :\ j \in J \setminus \{ b \}\, \}$, where
$\mathcal{H}^j = H^j \cap \mathcal{V}$.  A domain of the arrangement
$A$ is called {\it bounded in $\mathcal{V}$} if it is contained in a
suitable ball in $\mathcal{V}$.

Let $\Delta$ be a bounded domain and $\bar \Delta$ its closure.
Consider the subset $S \subset\bar \Delta$ of all maximally remote
points from the hyperplane $\mathcal{H}^{a}$. This subset is the
union of some open faces of $\bar \Delta$. The unique face $\Gamma
\subset S$ of highest dimension
 is called {\it the $\HH^a$-external supporting face} of $\Delta$.

Let $\Gamma$ be of dimension $m$, then there is a unique
$m$-dimensional edge $L$ of the projective arrangement $A$ which
contains $\Gamma$. The edge $L$ is called {\it the $\HH^a$-external
supporting edge of $\Delta$.} The triple $(L,H^a,H^b)$ is a
parallelism in $A$.

\begin{lemma}
\label{volume of parallelism} Let $(L,H^a,H^b)$ be a parallelism in
$A$. Then the number of bounded domains with $\HH^{a}$-external
supporting edge $L$ is equal to the discrete volume of the
parallelism $(L,H^{a},H^{b})$.
\end{lemma}

The lemma
is proved in Section \ref{Proof of
volume of parallelism}.

\section{Dual pairs} \label{dualarrangement}

\subsection{Admissible pairs} \label{admissible pairs}
Let $N$ be a natural number, $N \geq 3$. Let $k,n$ be natural
numbers such that $k+n+1=N$ and  $1 \leq k,n \leq N-2$.

Consider the vector space $ \R^{N+1}$ and its dual space. Let $\{
e_1 ,
 \ldots ,$ $ e_{N+1} \}$ be the standard basis of $\R^{N+1}$ and $\{
 e^1 , \ldots , e^{N+1} \}$ the dual basis of the dual space.
 Denote $\R^{N+1}$ by $\X$, denote the dual space by $\X^{\prime}$.
Denote $J = \{ 1 , \ldots , N+1 \}$.

Let ${\mathbb W} \subset \X$ be a vector subspace of dimension
$k+1$. Let $\W'\subset \X'$ be the annihilator of $\W$.  The
subspace $\W'$  is of dimension $n+1$.

The set of linear functions $\{\, e^j \vert_{\W}\ :\  j \in J\,\}$ spans
the dual space of $\W$.
Similarly, the set of linear functions
$\{\, e_j\vert_{\W'}\ :\  j \in J\,\}$ spans the dual
space of $\W'$.

Assume that for any $a,b \in J$, $a\neq b$, the functions $e^a
\vert_{\W}$ and $e^b \vert_{\W}$ are not proportional, and the
functions $e_a \vert_{\W'}$ and $e_b \vert_{\W'}$ are not
proportional.

The pair $\tau = (\X , \W )$ with this property will be called {\it
an admissible pair in $\X$}.  Similarly, the pair
$\tau' = (\X' , \W' )$ with this property will be called {\it an
admissible pair in $\X'$}. The pairs $\tau $ and
$\tau'$ will be called {\it dual}.

\subsection{The arrangements and matroid of an admissible pair}\label{duality}
Let $\tau = (\X , \W )$ be an admissible pair. For
$j \in J$, denote $E^j = \{ \,x \in \W\ :\ e^j (x) = 0\,\}$. These are
vector subspaces of $\W$ of codimension one.

Denote
\begin{enumerate}
\item [$\bullet$] $V = \Pee (\W)$, the projective space of dimension
$k$,
\item [$\bullet$]  $H^j = \Pee (E^j)$, $j \in J$, projective
hyperplanes in $V$,
\item [$\bullet$] $A [\tau] = \{\,H^j\ :\  j \in J \,\}$, the arrangement of
projective hyperplanes in $V$.
\end{enumerate}
Denote
\bea
\mathcal{J}\ =\ \{\,1, \ldots , N \,\}\ =\ J \setminus \{\,N+1\,\}\ .
\eea
For $j \in \mathcal{J}$, the rational function $f^j = {e^j}/
{e^{N+1}} $ restricted to $\W$ is regular on $\W \setminus E^{N+1}$
and homogeneous of degree zero. Thus, $f^j$ is a well-defined degree
one polynomial on the affine space
\bea
V \setminus H^{N+1}\ =\ \Pee (\W) \setminus \Pee (E^{N+1})\ .
\eea
Denote
\begin{enumerate}
\item [$\bullet$] $\mathcal{V} = V \setminus H^{N+1}$, the affine space of dimension
$k$,
\item [$\bullet$]  $\mathcal{H}^j = \{\,x \in \mathcal{V} \ :\ f^j (x) = 0\, \}$, $j \in \mathcal{J}$,
affine hyperplanes in $\mathcal{V}$,
\item [$\bullet$] $\mathcal{A}[\tau] = \{\,\mathcal{H}^j\ :\  j \in \mathcal{J}\, \}$,
the arrangement of
affine hyperplanes in $\mathcal{V}$.
\end{enumerate}
Observe that $\mathcal{H}^j = H^j \cap \mathcal{V}$, $j\in
\mathcal{J}$.

The set  $\mathcal F[\tau] = \{ \,f^j\ :\ j \in \mathcal{J}\, \}$ of
degree one polynomials on $\mathcal V$ will be called {\it the
arrangement of polynomials associated to $\tau$}.

For dual admissible pairs $\tau$ and $\tau'$, the corresponding pairs
of objects: \ $A [\tau]$ and $A [\tau']$,\  $\mathcal{A} [\tau]$ and
$\mathcal{A} [\tau']$,\  $\mathcal{F} [\tau]$ and $\mathcal{F}
[\tau']$\ -\ will be called {\it dual}.

Introduce {\it the matroid  of $\tau$}, denoted $M[\tau]$, as the
matroid of the collection of vectors $\{\, e^j \vert_{\W}\ :\ j \in  J
\,\}$ in the dual space of $\W$.

Observe that the matroid of $\tau$ is the matroid of the arrangement
$A[\tau]$.

\subsection{The value of a polynomial $f^j:\mathcal{V} \to \R$
 at a vertex} Let $P \in \mathcal{V}$ be a vertex of the arrangement
 $\mathcal{A}[\tau]$ and let $X$ be the flat associated to $P$
 in the matroid $M [\tau]$.
Let $Y=\{j_1,\dots,j_k\} \subset X$ be a maximal independent subset.
Let  $j\in \mathcal{J}$. Consider two vectors \bea v =(e^{j_1}\wedge
e^{j_2}\wedge \dots \wedge e^{j_k}\wedge e^{N+1})\vert_{\W} \ ,
\qquad v'=(e^{j_1}\wedge e^{j_2}\wedge \dots \wedge e^{j_k}\wedge
e^{j})\vert_{\W} \eea of the one-dimensional vector space
$\bigwedge^{k+1} \W^*$, where $\W^*$ is the dual space of $\W$. The
first vector is nonzero. Let $c\in\R$ be the coefficient of
proportionality: $v' = c v$.

\begin{lemma}
\label{value of polyn lemma}
The value of the polynomial $f^j$ at
$P$ is equal to $c$.
\hfill
$\square$
\end{lemma}

\subsection{A coordinate description}\label{coordinate description} Let $\tau = (\X , \W )$
and $\tau' = (\X^{\prime} , \W^{\prime} )$ be dual admissible pairs.

Let $w_1 , \ldots , w_k , w_{k+1} , w_{k+2} , \ldots , w_{N+1} $ be
any basis of $\X$ such that $w_1 , \ldots , w_k , w_{k+1}$ is a
basis of $\W$. Consider the dual basis $w^1 , \ldots , w^k , w^{k+1}
, w^{k+2} , \ldots , w^{N+1}$ of $\X^{\prime}$. Then $w^{k+2} ,
\ldots , w^{N+1}$ is a basis of $\W^{\prime}$.

Let $w_i = \sum_{l=1}^{N+1}b_i^l e_l$ and  $w^i =
\sum_{l=1}^{N+1}c_l^i e^l$, for  $i=1,2, \dots , N+1$. Denote
$\mathfrak{B} = (b_i^l)$, $\mathfrak{C} = (c^i_l)$. We have
$\mathfrak{C} = (\mathfrak{B}^T)^{-1}$.

Introduce the $(k+1) \times (N+1)$-matrix $B$ and the $(n+1) \times
(N+1)$-matrix $C$ by
$$
B = \left(
\begin{array}{lllllll}
b_1^1 & b_1^2 & \cdots & b_1^j & \cdots & b_1^N & b_1^{N+1}\\
\noalign{\medskip}
b_2^1 & b_2^2 & \cdots & b_2^j & \cdots & b_2^N & b_2^{N+1}\\
\noalign{\medskip}
\vdots & \vdots & & \vdots & & \vdots & \vdots \\
\noalign{\medskip}
b_k^1 & b_k^2 & \cdots & b_k^j & \cdots & b_k^N & b_k^{N+1} \\
\noalign{\medskip}
b_{k+1}^1 & b_{k+1}^2 & \cdots & b_{k+1}^j & \cdots & b_{k+1}^N & b_{k+1}^{N+1}\\
\end{array}
\right)
$$
and
$$
C = \left(
\begin{array}{lllllll}
c_1^{k+2} & c_2^{k+2} & \cdots & c_j^{k+2} & \cdots & c_N^{k+2} & c_{N+1}^{k+2} \\
\noalign{\medskip}
c_1^{k+3} & c_2^{k+3} & \cdots & c_j^{k+3} & \cdots & c_N^{k+3} & c_{N+1}^{k+3} \\
\noalign{\medskip}
\vdots & \vdots & & \vdots & & \vdots & \vdots \\
\noalign{\medskip}
c_1^{N} & c_2^{N} & \cdots & c_j^{N} & \cdots & c_N^{N} & c_{N+1}^{N}\\
\noalign{\medskip} c_1^{N+1} & c_2^{N+1} & \cdots & c_j^{N+1} &
\cdots & c_N^{N+1} & c_{N+1}^{N+1}
\end{array}
\right) .
$$
The matrices $B$ and $C$ are parts of the matrices $\mathfrak{B}$
and $\mathfrak{C}$, respectively. Clearly, $\hbox{rank } B = k+1$
and $\hbox{rank } C = n+1$.

Let $x^1 , \ldots , x^{k+1} : \W \longrightarrow \R$ be the
coordinate functions with respect to the basis $w_1 , \ldots ,
w_{k+1}$, $x^i (w_j) = \delta^i_j , \quad i,j = 1, \ldots , k+1$.
Observe that for any $j=1, \dots, k+1$, we have $x^j =
w^j\vert_{\W}$. Then for any $j\in J$, we have
$$
e^j\vert_{\W} = b^j_1 x^1 + \dots + b^j_{k+1} x^{k+1} \ .
$$
Similarly, let $x_1 , \ldots , x_{n+1} : \W^{\prime} \longrightarrow
\R$ be the coordinate functions
 with respect to the basis $w^{k+2} , \ldots , w^{N+1}$,
$x_i (w^{k+j+1}) = \delta_i^j , \quad i,j = 1, \ldots , n+1$.
Observe that for any $j=1, \dots, n+1$, we have $x_j =
w_{k+j+1}\vert_{\W}$. Then for any $j\in J$, we have
$$
e_j\vert_{\W'} = c_j^{k+2} x_1 + \dots + c_j^{N+1} x_{n+1} \ .
$$
Thus, the columns of $B$ and $C$ describe the coordinates of the
functions $e^j\vert_{\W}$ and $e_j\vert_{\W'}$.

Denote by
 $B [j_1 ,j_2, \ldots , j_{k+1}]$, the determinant of the
$(k+1) \times (k+1)$-submatrix of $B$ formed by the rows $1,2,
\ldots ,k+1$ and columns $j_1 < j_2 < \cdots < j_{k+1}$. Denote by
 $C [l_1 ,l_2, \ldots , l_{n+1}]$, the determinant of the
$(n+1) \times (n+1)$-submatrix of $C$ formed by the rows $1,2,
\ldots ,n+1$ and columns $l_1 < l_2 < \cdots < l_{n+1}$.

The values of functions $f^j$, $j \in \mathcal{J}$, at the vertices
of the affine arrangement $\mathcal{A}[\tau]$ in $\mathcal{V}$ can
be computed in terms of the minors of the matrix $B$ as follows.

\begin{lemma} 
\label{value of fj}
Let $P = \mathcal{H}^{j_1} \cap \cdots \cap \mathcal{H}^{j_k}$ be a
vertex of the arrangement $\mathcal{A}[\tau]$ lying in
$\mathcal{V}$, $j_1 < \cdots < j_k$. Let $j \in \mathcal{J}
\setminus \{j_1 , \ldots , j_k \}$. Assume that $j_1< \dots < j_m <
j< j_{m+1} < \dots < j_k$ where $m \leq k$.  Then \bea
\phantom{aaaaaaaaaaaaaaaaa} f^j (P) = (-1)^{k + m  }  \, \frac{B
[j_1 , \dots , j_m, j , j_{m+1}, \ldots , j_k]}{B [j_1 , \dots , j_k
, N+1]} \ . \phantom{aaaaaaaaaaa}
 \square
\eea
\end{lemma}

A similar statement holds for the arrangement $\mathcal{A}[\tau']$
and matrix $C$.

\subsection{Relation between minors of $B$ and $C$}

\begin{theorem} [\cite{M}, pp. 165--169]
\label{minors} ${}$  \ Let $1 \leq j_1 < \cdots < j_{k+1} \leq N+1$
be a subset and let $1 \leq j_{k+2} < \cdots < j_{N+1} \leq N+1$ be
the complementary subset. Then
$$
B [j_1 , \ldots , j_{k+1}] = (-1)^\sigma \cdot \det \mathfrak{B}
\cdot C [j_{k+2} , \ldots , j_{N+1}]\ ,
$$
where
$$
\sigma = 1+2+ \cdots + (k+1) + j_1 + j_2 + \cdots + j_{k+1}\ .
$$
\end{theorem}

\subsection{Dual matroids}
Let $M=(J, \mathcal{I}\,)$ be a matroid. A subset $X \subset J$
belongs to $\mathcal{I}$ if and only if $X$ is contained in at least
one basis of $M$. In this way, a matroid is characterized by its
collection of bases.

The {\it dual} of $M$ is the matroid $M'$ on the same ground set $J$,
whose bases are the complements in $J$ of the bases of $M$. Evidently,
$(M')' = M$.

\begin{theorem} [\cite{O}, \cite{B}] \label{mmstar}
${}$ \ Let $M$ and $M'$ be dual matroids. Then
\begin{enumerate}
\item[$\bullet$]
For any subset $X \subset J$, $\text{\rm corank}_M\,  X \ =\
\text{\rm nullity}_{M'} \, \hat{X}$,
\item[$\bullet$]
 ${\rm T} (M; x, y) = {\rm T} (M'; y, x)$.
\end{enumerate}
\end{theorem}

\begin{lemma} \label{b10cor}
Let $M$ and $M'$ be dual matroids, $|J|\geq 2$. Then $b^{10}_M =
b^{10}_{M'}$. \hfill $\square$
\end{lemma}
The lemma follows from Theorems \ref{b10} and \ref{mmstar}.

\begin{lemma} \label{delcont}
Let $M$ and $M'$ be dual matroids. Let $X \subset J$, $|X|<|J|$.
Then the matroids $M / X$ and $M' - X$ are dual.
\end{lemma}
The lemma
is proved in Section \ref{proof of delcont}.

\begin{theorem} \label{volumelemma}
Let $X$ be a spacious flat in $M$, $1<|X|< |J|-1$. Then the
complement of $X$, the set  $\hat{X} = J \setminus X$, is a
spacious flat in $M'$. Furthermore,
\begin{enumerate}
\item[$\bullet$] ${\rm l}_{M} X = {\rm w}_{M'} \hat{X}$,
\item[$\bullet$] ${\rm w}_{M} X = {\rm l}_{M'} \hat{X}$,
\item[$\bullet$] ${\rm vol}_{M} X = {\rm vol}_{M'} \hat{X}$.
\end{enumerate}
\end{theorem}
Theorem \ref{volumelemma} is proved in Section \ref{proof of
volumelemma}.

Recall that for  a parallelism $(X , a , b)$  in $M$, we denote by
$\hat X(a,b)$ the complement of $X \cup \{a,b \}$ in $J$.

\begin{theorem} \label{triplevol}
Let $(X , a , b)$ be a parallelism in $M$, $|X|<|J|-2$. Assume that
${\rm corank}_M  (J \setminus \{ a,b \}) =0$. Then $(\hat X(a,b) , a
, b )$ is a parallelism in $M'$. Furthermore,
\begin{enumerate}
\item[$\bullet$] ${\rm l}_{M} X = {\rm w}_{M'} (\hat X(a,b), a , b)$,
\item[$\bullet$] ${\rm w}_{M} (X , a , b) = {\rm l}_{M^{\prime}} \hat X(a,b)$,
\item[$\bullet$] ${\rm vol}_{M} (X , a , b) =
{\rm vol}_{M^{\prime}} (\hat X(a,b) , a , b)$.
\end{enumerate}
\end{theorem}
Theorem \ref{triplevol} is proved in Section \ref{proof of triplevol}.

\begin{lemma} \label{arrangement matroids dual}
Let $\tau$ and $\tau'$ be dual admissible pairs. Then the matroids
$M [\tau]$ and $M [\tau']$ are dual.  \hfill $\square$
\end{lemma}
The lemma is a corollary of  Theorem \ref{minors}.

\subsection{Bounded domains, edges, and parallelisms of dual
arrangements}\label{domains edges parallelisms}

\begin{lemma}
The number of bounded domains of the arrangements $A[\tau]$ and $A
[\tau']$ are equal. \hfill $\square$
\end{lemma}
The lemma follows from Theorem \ref{infinity} and Lemmas
\ref{b10cor}, \ref{arrangement matroids dual}.

For dual admissible pairs $\tau$ and $\tau'$, write
\begin{enumerate}
\item[$\bullet$] $A [\tau] = \{H^j : j \in J \}$,\qquad {}
 $A [\tau'] = \{ H_j : j \in J\}$,
\item[$\bullet$] $\mathcal{A} [\tau] = \{ \mathcal{H}^j : j \in
\mathcal{J}\}$,\qquad $\mathcal{A} [\tau'] = \{ \mathcal{H}_j : j
\in \mathcal{J}\}$,
\item[$\bullet$] $\mathcal{F}[\tau] = \{ f^j : j \in \mathcal{J}\}$,
\qquad $\mathcal{F}[\tau'] = \{ f_j : j \in \mathcal{J}\}$.
\end{enumerate}

Let $L$ be a spacious edge of the arrangement $A[\tau]$ and $X$ the
associated flat in the matroid $M[\tau]$. Denote by $\hat L$ the
edge $ \cap_{j \in \hat X} H_j$ of the arrangement $A[\tau']$.

\begin{lemma} \label{dual edges}
Let $L$ be a spacious edge of $A[\tau]$. Assume that $L$ is not a
hyperplane. Then $\hat L$ is a spacious edge in $A[\tau']$.
Furthermore,
\begin{enumerate}
\item[$\bullet$] ${\rm l}_{A[\tau]} L = {\rm w}_{A [\tau']} \hat{L}$,
\item[$\bullet$] ${\rm w}_{A[\tau]} L = {\rm l}_{A [\tau']} \hat{L}$,
\item[$\bullet$] ${\rm vol}_{A[\tau]} L = {\rm vol}_{A [\tau']} \hat{L}$. \hfill
$\square$
\end{enumerate}
\end{lemma}
The lemma follows from Theorem \ref{volumelemma}.

The edges $L$ and $\hat L$ will be called {\it dual}.

Let $(L , H^{a} , H^{b})$ be a parallelism in the arrangement
$A[\tau]$. Let $X$ be the flat in the matroid $M[\tau]$ associated
to the edge $L$. Denote by $\hat L(H^a,H^b)$ the edge $\cap_{j \in
\hat{X}(a,b)} H_j \,$ in the arrangement $A[\tau']$.

\begin{lemma} \label{dual parallelisms}
Let $(L , H^{a} , H^{b})$ be a parallelism in $A[\tau]$. Then $(\hat
L(H^a,H^b) , H_{a} , H_{b})$ is a parallelism in $A[\tau']$.
Furthermore,
\begin{enumerate}
\item[$\bullet$] ${\rm l}_{A[\tau]} L  = {\rm w}_{A [\tau']}
(\hat L(H^a,H^b) , H_{a} , H_{b})$,
\item[$\bullet$] ${\rm w}_{A[\tau]} (L , H^{a} , H^{b}) = {\rm l}_{A
[\tau']} \hat L(H^a,H^b)$,
\item[$\bullet$] ${\rm vol}_{A[\tau]} (L , H^{a} , H^{b}) = {\rm
vol}_{A [\tau']} (\hat L(H^a,H^b) , H_{a} , H_{b})$. \hfill
$\square$
\end{enumerate}
\end{lemma}
The lemma follows from Theorem \ref{triplevol}.

The parallelisms $(L , H^{a} , H^{b})$ and $( \hat L(H^a,H^b) ,
H_{a} , H_{b})$ will be called {\it dual}.

\subsection{Relation between the values of $f^j$ and $f_j$}

\begin{lemma}
\label{product equals negative one} Let $(L , H^j , H^{N+1})$ be a
parallelism in $A[\tau]$ and
\linebreak
$(\hat{L} (H^j , H^{N+1}) , H_j ,
H_{N+1})$ the dual parallelism in $A[\tau']$. The product of the
value of $f^j$ on $L \setminus H^{N+1}$ and the value of $f_j$ on
$\hat{L} (H^j , H^{N+1}) \setminus H_{N+1}$ is $-1$.
\end{lemma}
Lemma \ref{product equals negative one} is proved in Section
\ref{proof of product equals negative one}.

\section{Determinant formula}\label{Sec DET}
\subsection{Weighted arrangements}
Let $A = \{H^j : j \in J \}$ be a projective arrangement. A set of
numbers $\alpha = \{ \alpha_j : j \in J\}$ with the property
$\sum_{j \in J} \alpha_j = 0$ will be called
 {\it a set of weights}, where $\alpha_j$ is {\it the weight of} $H^j$. The pair
$(A,\alpha)$ will be called {\it a weighted arrangement}.

{\it The weight of an edge} $L$ of $A$ is the sum, $\alpha (L)$, of
the weights of the hyperplanes that contain $L$.

Let $A [\tau] = \{ H^j : j \in J \}$ and $A [\tau^{\prime}] = \{ H_j
: j \in J \}$ be dual arrangements.  Let $\alpha = \{ \alpha_j : j
\in J\}$ be a set of numbers with the property $\sum_{j \in J}
\alpha_j = 0$. Then  $( A [\tau] , \alpha )$ and $( A
[\tau^{\prime}] , \alpha )$ will be called {\it dual weighted
arrangements}.

Let $( A [\tau] , \alpha )$ be a weighted arrangement. For $j \in
\mathcal{J}$, the number $\alpha_j$ is also called the {\it weight
of the affine hyperplane} $\mathcal{H}^j$. The associated affine
arrangement $\mathcal{A} [\tau]$ with weight $\alpha_j$ assigned to
the hyperplane $\mathcal{H}^j$ for any $j \in \mathcal{J}$, will be
called  {\it the associated weighted affine arrangement} and denoted
$(\mathcal{A} [\tau] , \alpha)$.

We always assume that for any $j \in \mathcal{J}$,
 the weight $\alpha_j$ of the hyperplane $\mathcal{H}^j$ is a positive real
number.

\subsection{Period matrix of a weighted arrangement}
Let $( \mathcal{A} [ \tau ] , \alpha )$ be a weighted affine
arrangement. We have $\mathcal{A} [\tau] = \{ \mathcal{H}^j : j \in
\mathcal{J}\}$, where $\mathcal J = \{1,\dots, N\}$. Consider
$\mathcal J$ as an ordered set with the standard order.

In \cite{DT}, for an ordered affine arrangement,  a collection $
\beta${\bf kbc}($\mathcal{A} [\tau]$) of ordered $k$-tuples $B =
(\mathcal{H}^{j_1} , \mathcal{H}^{j_2} , \ldots ,
\mathcal{H}^{j_k})$ is defined. For any tuple of that collection,
the intersection $\cap_{l=1}^k \mathcal{H}^{j_l}$ is a vertex. The
number of tuples in that collection is equal to the number $\beta =
\beta (\mathcal A[\tau])$ of bounded domains of the affine
arrangement.

Elements of the collection are called $\beta${\bf kbc}-bases. The
collection itself   is ordered lexicographically, $  \beta${\bf
kbc}($\mathcal{A} [\tau]$) $ = \{ B_1,\dots, B_\beta\}$.

\subsubsection{Logarithmic $k$-forms} \label{enumkforms}
For $B = (\mathcal{H}^{j_1} , \mathcal{H}^{j_2} , \ldots ,
\mathcal{H}^{j_k}) \in \beta \hbox{\bf kbc}(\mathcal{A} [\tau])$,
define the associated flag of edges \bea \xi (B) = ( L_B^0 \subset
L_B^1 \subset \cdots \subset L_B^k = \mathcal{V} )\ , \eea where
$L_B^i = \mathcal{H}^{j_{i+1}} \cap \mathcal{H}^{j_{i+2}} \cap
\cdots \cap \mathcal{H}^{j_k}$ and  $\dim L_B^i = i$.

Let $\mathcal{F} [\tau] = \{f^j : j \in \mathcal{J} \}$ be the
polynomial arrangement associated to $\tau$. Let $L\subset
\mathcal{V}$ be an edge of $\mathcal{A} [\tau]$. Assign to $L$ the
differential $1$-form \bea \omega (L) = \sum_{j, \ L \subset
\mathcal H^j} \alpha_j \frac{d f^j}{f^j}\ . \eea Assign to every $B =
(\mathcal{H}^{j_1} , \mathcal{H}^{j_2} , \ldots ,
\mathcal{H}^{j_k})$
 the differential $k$-form
$$
\phi (B) = \omega (L_{B}^0 ) \wedge \omega (L_{B}^1 ) \wedge \cdots
\wedge \omega (L_{B}^{k-1} )\  .
$$
Thus, we get an ordered set of differential $k$-forms
$$
\Psi = \{ \phi^1 , \phi^2 , \ldots , \phi^{\beta} \}\  ,
$$
where $\phi^i$ is the form corresponding to the $i$-th element of
the collection $  \beta${\bf kbc}($\mathcal{A} [\tau]$).

\subsubsection{The  $\beta${\bf kbc}-enumeration of bounded domains}
\label{enumdomains} Let $\xi = ( L^0 \subset L^1 \subset \cdots
\subset L^k )$ be a flag of edges of $\mathcal{A} [\tau]$ with $\dim
L^j = j$ for all $j$. Let $\Delta$ be a bounded domain of
$\mathcal{A} [\tau]$ and $\overline \Delta$ its closure.
 The flag $\xi$ is said
to be {\it adjacent to} $\Delta$, if $\dim\, ( L^j \cap
\overline{\Delta} \, ) = j$, for all $j$.

Denote by $\text{\tt Ch} (\mathcal{A} [\tau])$ the set of bounded
domains of the arrangement $\mathcal{A} [\tau]$. In \cite{DT},
\linebreak
 a bijection
$$
C : \beta \text{\bf kbc}(\mathcal{A} [\tau] ) \longrightarrow
\text{\tt Ch} (\mathcal{A} [\tau])
$$
is defined such that for any $B \in \beta \,\text{\bf
kbc}(\mathcal{A} [\tau] )$, the associated flag $\xi (B)$ is
adjacent to the bounded domain $C (B)$.

Thus, one has $\text{\tt Ch} (\mathcal{A} [\tau]) = \{ \Delta_1 ,
\Delta_2 , \ldots , \Delta_{\beta} \}$, where
$$\Delta_s = C (B_s)\ , \quad s= 1 , 2 , \ldots , \beta \ .$$ This is called the
$\beta \text{\bf kbc}(\mathcal{A} [\tau] )$-{\it ordering} of the
bounded domains of $\mathcal{A} [\tau]$.

\subsubsection{Orientation of bounded domains} \label{orient}
Let $\Delta = C (B)$ and $\xi (B) = ( L^0 \subset L^1 \subset \cdots
\subset L^k )$. The flag $\xi (B)$ is adjacent to the domain $\Delta
= C(B)$ and defines the {\it intrinsic orientation} of $\Delta$, see
Section 6.2 in \cite{V2}.

The intrinsic orientation is the orientation of the unique
orthonormal frame \linebreak $\{ v_1, v_2, \ldots , v_k \}$ where
$v_i$ is the unit vector originating from the vertex $L^0$
 in the direction of $L^i \cap \overline{\Delta}$.

\subsubsection{The period matrix} \label{the period matrix}
Let $( \mathcal{A} [\tau] , \alpha )$ be a weighted affine
arrangement and \linebreak $\mathcal{F}[\tau] = \{ f^j : j \in
\mathcal{J}\}$ the associated polynomial arrangement. Consider the
multi-valued function
$$
U^{\alpha}\ = \ \prod_{j \in \mathcal{J}}\ (f^j)^{\alpha_j}\  :\
 \mathcal{V} \longrightarrow \C\ .
$$
Fix a uni-valued branch of each $(f^j)^{\alpha_j}$ on each bounded
domain $\Delta$. The $\beta \times \beta $-matrix
$$
 \text{PM} \,(\mathcal{A} [\tau] , \alpha ) \ =\  \left(\, \int_{\Delta_s}
U^{\alpha} \, \phi^t \,\right)
$$
is called {\it the period matrix} of the weighted affine arrangement
$( \mathcal{A} [\tau] , \alpha)$. In this matrix the differential
forms are ordered as in Section \ref{enumkforms}, the domains are
ordered as in Section \ref{enumdomains}, the domains are oriented as
in Section \ref{orient}, and in each integral of the matrix the
uni-valued branch of $U^\alpha$ is chosen since the branches of each
of its factors were chosen.

Denote the determinant of the period matrix by $ \text{D}
(\mathcal{A} [\tau] , \alpha )$.

\medskip

\noindent
{\bf Remark.} The ordered set of differential $k$-forms described in
Section \ref{enumkforms},  the order on the set of bounded domains
as in Section \ref{enumdomains}, the orientation on bounded domains
as in Section \ref{orient} will be called {\it canonical}.

\subsection{Determinant of the period matrix}
\subsubsection{Beta function}
\label{betafunction} Let $\mathcal{L}_- [\tau]$ be the set of all
edges of $A [\tau]$ lying in $H^{N+1}$ and $\mathcal{L}_+ [\tau]$
the set of all other edges. The {\it beta function of the weighted
affine arrangement} $(\mathcal{A} [\tau] , \alpha )$ is defined in
\cite {V1} as
$$
B (\mathcal{A} [\tau] , \alpha ) \ =\ \frac{\prod_{L \in
\mathcal{L}_+ [\tau]} \Gamma (\alpha (L) + 1 )^{{\rm vol}_{A[\tau]}
(L)}}{\prod_{L \in \mathcal{L}_- [\tau]} \Gamma (- \alpha (L) + 1
)^{{\rm vol}_{A[\tau]} (L)}} \ {}\ ,
$$
where $\Gamma$ is Euler's gamma function.

\subsubsection{Critical values}
\label{criticalvalues} Let $\Delta$ be a bounded domain of
$\mathcal{A} [\tau]$. Let $j \in \mathcal J$. Let $\Sigma$ be the
$\mathcal{H}^j$-external supporting face of $\Delta$. The value of
the chosen
 branch of ${(f^j)}^{\alpha_j}$ on $\Sigma$
is called {\it the critical value of $(f^j)^{\alpha_j}$ on $\Delta$}
and denoted by $ c \,((f^j)^{\alpha_j} , \Delta )$.

\subsubsection{Evaluation of the determinant}
\begin{theorem} [\cite{V1}, \cite{DT}]
\label{thm evaluation} The determinant of the period matrix is given
by the following formula:
$$
{\rm D} (\mathcal{A} [\tau] , \alpha ) = B (\mathcal{A} [\tau],
\alpha ) \ \cdot \prod_{\stackrel{j \in \mathcal{J}}{\Delta \in \,
{\tt Ch \mathcal{A} [\tau]}}} c \,({(f^j)}^{\alpha_j} , \Delta ) \
.
$$
\end{theorem}

\subsection{Special choice of branches}
The definition of the period matrix ${\rm PM} (\mathcal{A} [\tau] ,
\alpha )$ involves the choice of branches of each $(f^j)^{\alpha_j}$
on each bounded domain of $\mathcal{A}[\tau]$. In this paper, we
choose the branches as follows.

Let $(L , H^j , H^{N+1})$ be a parallelism in $A [\tau]$. On all
bounded domains of $\mathcal{A} [\tau]$ whose
$\mathcal{H}^j$-external supporting edge is $L$, choose the argument
of $f^j$ to be the same.

Then for any two bounded domains $\Delta$ and $\Delta'$ with the
same $\mathcal{H}^j$-external supporting edge $L$, the values $c
\,({(f^j)}^{\alpha_j} , \Delta)$ and $c \,({(f^j)}^{\alpha_j} ,
\Delta')$ will be equal. Denote this common value by $c
\,({(f^j)}^{\alpha_j} , L )$.

Repeating this process for all parallelisms $(L , H^j , H^{N+1})$ in
$A[\tau]$ gives a choice of branches of each $(f^j)^{\alpha_j}$ on
each bounded domain of $\mathcal{A}[\tau]$.

Such a choice of branches will be called {\it special}.

\begin{lemma} For a special choice of branches, construct the period matrix
\linebreak
${\rm PM} (\mathcal{A} [\tau] , \alpha )$. Then its determinant is given
by the formula:
$$
{\rm D} (\mathcal{A} [\tau] , \alpha ) = B (\mathcal{A} [\tau] ,
\alpha ) \cdot \prod c \,({(f^j)}^{\alpha_j} , L )^{{\rm
vol}_{A[\tau]} (L , H^j , H^{N+1})} \ ,
$$
where the product is taken over all parallelisms $(L , H^j ,
H^{N+1})$. \hfill $\square$
\end{lemma}

The lemma follows from Lemma \ref{volume of parallelism} and Theorem
\ref{thm evaluation}.

\subsection{Associated period matrices of dual arrangements}\label{associated period matrices}

Let $(\mathcal{A} [\tau] , \alpha )$ and $(\mathcal{A} [\tau'] ,
\alpha )$ be dual weighted affine arrangements. For each of them we
can define period matrices and calculate their determinants. The
period matrices depend on the choice of branches of functions
$(f^j)^{\alpha_j}$ and $(f_j)^{\alpha_j}$ on bounded domains of
those arrangements.

In this section we define the associated choices of branches in such
a way that the determinants of period matrices will be related.

Consider $\mathcal{J} = \{1,\dots,N\}$ with the standard order. Then
we have the canonical ordered set of differential $k$-forms
associated with $\mathcal A[\tau]$, ordering on the set of bounded
domains of $\mathcal{A} [\tau]$, and orientation on each bounded
domain. We also have the canonical ordered set of differential
$n$-forms associated with $\mathcal A[\tau']$, ordering on the set
of bounded domains of $\mathcal{A} [\tau']$, and orientation on each
bounded domain.

Take an arbitrary  special choice of branches of the functions
$(f^j)^{\alpha_j}$ on bounded domains of the arrangement $\mathcal
A[\tau]$. We will define now the associated special choice of
branches of functions $(f_j)^{\alpha_j}$ on bounded domains of the
arrangement $\mathcal A[\tau']$.

Let $(L , H^j , H^{N+1})$ be a parallelism in $A[\tau]$ and
$(\hat{L} (H^j , H^{N+1}) , H_j , H_{N+1})$ the dual parallelism in
$A[\tau']$. Suppose that in the definition of ${\rm PM} (\mathcal{A}
[\tau] , \alpha )$, $\theta$ is the chosen argument of $f^j$ on a
bounded domain with $\mathcal{H}^j$-external supporting edge $L$. On
each bounded domain of $\mathcal{A} [\tau^{\prime}]$ with
$\mathcal{H}_j$-external supporting edge $\hat{L} (H^j , H^{N+1})$,
choose the argument of $f_j$ to be $ - \theta +  \pi $, see Lemma
\ref{product equals negative one}.

With this choice of branches, define the period matrix, ${\rm PM}
(\mathcal{A} [\tau^{\prime}] , \alpha )$, of $(\mathcal{A}
[\tau^{\prime}] , \alpha )$.

The period matrices ${\rm PM} (\mathcal{A} [\tau] , \alpha )$ and
${\rm PM} (\mathcal{A} [\tau^{\prime}] , \alpha )$ will be called
{\it associated}.

The main result of this paper is the following theorem.

\begin{theorem}\label{main theorem}
The product of determinants of the associated period matrices is
given by the formula:
$$
{\rm D} ( \mathcal{A} [\tau] , \alpha ) \cdot {\rm D} ( \mathcal{A}
[\tau^{\prime}] , \alpha ) = \left[ \frac{\prod_{j \in \mathcal{J} }
e^{\pi i \alpha_j} \Gamma (\alpha_j + 1)}{\Gamma ( \sum_{ j \in
\mathcal{J}} \alpha_j + 1 )} \right]^{\beta} \ ,
$$
where $\beta$ is the number of bounded domains in $\mathcal{A}
[\tau]$.
\end{theorem}

Recall that the number of bounded domains in $\mathcal{A} [\tau]$ is
equal to the number of bounded domains in $\mathcal{A} [\tau']$.

The theorem is proved in Section \ref{proof of main theorem}.

\section{proofs}
 \label{proofs}
\subsection{Proof of Lemma \ref{volume of parallelism}}
\label{Proof of volume of parallelism}
Let $(L, H^a , H^b)$ be a parallelism of the arrangement $A$ in a
projective space $V$. Let $\Sigma$ be a domain of the induced
arrangement $A_L$ on $L$, bounded with respect to the hyperplane $H^b \cap
L$. It is enough to show that the number of bounded domains with
$\HH^a$-external supporting face $\Sigma$ is ${\rm w}_A (L, H^a ,
H^b)$, the discrete width of the parallelism $(L, H^a , H^b)$.

Denote by $X$ the flat associated to $L$ in the matroid, $M$, of the
arrangement $A$. By definition, ${\rm w}_A (L, H^a , H^b) =
b^{10}_{M - \hat{X} (a,b)}$.

Suppose that $L$ is a vertex. Then, by Theorem \ref{infinity},
$b^{10}_{M - \hat{X} (a,b)}$ is the number of bounded domains formed
by the set of hyperplanes
\bea
A^{(L, H^a, H^b)} = \{ H^j : j \in X \cup \{ a , b\}\}\ .
\eea
Clearly, bounded domains of the arrangement $\mathcal A$  with
$\HH^a$-external supporting face $\Sigma$ are in one-to-one
correspondence with bounded domains formed by $A^{(L, H^a, H^b)}$.
Thus, the number of bounded domains with $\HH^a$-external supporting
face $\Sigma$ is
\linebreak
 ${\rm w}_A (L, H^a , H^b)$.

If $L$ is not a vertex, consider a subspace $V' \subset V$, $\dim V'
= {\rm codim}\, L$, such that $V'$ intersects $L$ transversally.
Then $V'$ also intersects $H^a$ and $H^b$ transversally. Consider
the induced arrangement $\{H^j \cap V': j \in X \cup \{a,b\} \}$ on
$V'$. Each bounded domain of $A$ with $\HH^a$-external support
$\Sigma$ determines a unique bounded domain (with respect to the
hyperplane $H^b \cap V'$) of the new arrangement. The number of
which is $b^{10}_{M - \hat{X} (a,b)}$ by Theorem \ref{infinity}.
 \hfill $\square$

\subsection{Proof of Lemma \ref{delcont}} \label{proof of delcont}
Let $B \subset \hat{X}$ be a basis of $M / X$. Then there exists a
maximal independent subset $Y$ of $X$ such that $B \cup Y$ is a
basis of $M$.
We want to show that $\hat{X} \setminus B$ is a basis of $M' - X$.
That is, $\hat{X} \setminus B$ is a maximal independent subset of
$\hat{X}$ in $M'$.

Since $B \cup Y$ is a basis of $M$, $(\hat{X} \setminus B)\cup (X
\setminus Y)$ is a basis of $M'$. Hence, $\hat{X} \setminus B$ is
independent.
Now, $|B|={\rm corank}_M X = {\rm nullity}_{M'} \hat{X} = |\hat{X}|
- {\rm rank}_{M'} \hat{X}$. Therefore, $|\hat{X} \setminus B|={\rm
rank}_{M'} \hat{X}$. This shows that $\hat{X} \setminus B$ is
maximally independent in $\hat{X}$ in $M'$.

Conversely, let $B$ be a basis of $M^{\prime} - X$. Then $B$ is a
maximal independent subset of $\hat{X}$ in $M^{\prime}$. Hence there
exists an independent subset $Y$ of $X$ in $M'$ such that $B \cup Y$
is a basis of $M'$. Thus, $(\hat{X} \setminus B) \cup (X \setminus
Y)$ is a basis of $M$. Hence, $X \setminus Y$ is independent in $M$.
It remains to show that it is maximal independent.

Since $B \cup Y$ is a basis of $M'$ and $B$ is a maximal independent
subset of $\hat{X}$ in $M'$,
$
|Y|= {\rm corank}_{M'} \hat{X} = {\rm nullity}_M X
$.
Hence, $ |X \setminus Y|=|X|- {\rm nullity}_M X = {\rm rank}_M X
$. \hfill $\square$

\subsection{Flats of the dual matroid}
\begin{lemma} \label{dual flats}
Let $M = (J, \mathcal{I}\,)$ be a matroid. Let $X \subset J$ be a
subset, $1 < |X| < |J|-1$. Then $\hat{X}$ is a flat in the dual
matroid $M'$ if and only if the deletion $M - \hat{X}$ does not have
an isthmus.
\end{lemma}
\begin{proof}
Suppose that $M-\hat{X}$ does not have an isthmus. Then for every
$e \in X$, ${\rm rank}_{M - \hat{X}} X \setminus \{e\} = {\rm
rank}_{M-\hat{X}} X$. Hence, for every $e \in X$, ${\rm rank}_M X
\setminus \{e\} = {\rm rank}_M X$. That is, for every $e \in X$,
${\rm nullity}_{M} X \setminus \{ e\} = {\rm nullity}_{M} X -1$.
Thus, for every $e \in X$, ${\rm corank}_{M'} \hat{X} \cup \{ e\} =
{\rm corank}_{M'} \hat{X} -1$. Then for every $e \in X$, ${\rm
rank}_{M'} \hat{X} \cup \{ e\} = {\rm rank}_{M'} \hat{X} +1$. This
shows that the subset $\hat{X} \subset J$ is a flat in $M'$.

The converse follows by tracing the arguments backward.
\end{proof}

\subsection{Proof of Theorem \ref{volumelemma}} \label{proof of volumelemma}
\begin{lemma}\label{isthmus 1}
Let $X$ be a spacious flat in $M$, $1<|X|<|J|-1$. Then the deletion
$M - \hat{X}$ does not have an isthmus.
\end{lemma}
\begin{proof}
Let $e$ be an isthmus in the matroid $M - \hat{X}$. By Theorem
\ref{tuttepolynomial}, $ {\rm T} (M - \hat{X} ; x , y) = x \, {\rm
T} ((M - \hat{X}) / \{e\};x,y  )$. Hence, $ {\rm w}_M (X) =
b^{10}_{M - \hat{X}} = b^{00}_{(M - \hat{X}) / \{e\} }.$\ By
Theorem \ref{b10}, $ b^{00}_{(M - \hat{X}) /\{ e\} } = 0$ (since
$|X| \geq 2$). This contradicts our assumption that the flat $X$ is
spacious. Hence $M - \hat{X}$ does not have an isthmus.
\end{proof}

\begin{proof}[Proof of Theorem \ref{volumelemma}]
Let $M$ and $M^{\prime}$ be dual matroids. Let $X$ be a spacious
flat in $M$, $1<|X|<|J|-1$. It follows from Lemmas \ref{dual flats},
\ref{isthmus 1} that $\hat{X}$ is a flat in $M'$.

Furthermore, ${\rm l}_M X = b^{10}_{M/X} = b^{10}_{M' - X} = {\rm
w}_{M'} \hat{X}$ by Lemmas \ref{b10cor}, \ref{delcont}. Hence, ${\rm
vol}_M X = {\rm vol}_{M'} \hat{X}$.
\end{proof}

\subsection{Proof of Theorem \ref{triplevol}} \label{proof of triplevol}
\begin{lemma}\label{isthmus 2}
Let $(X,a,b)$ be a parallelism in $M$, $|X|<|J|-2$. Then the
deletion $M - \hat{X}(a,b)$ does not have an isthmus. \hfill
$\square$
\end{lemma}
The proof is similar to the proof of Lemma \ref{isthmus 1}.

\begin{proof}[Proof of Theorem \ref{triplevol}]
Let $(X,a,b)$ be a parallelism in $M$, $|X|<|J|-2$, \linebreak ${\rm
vol}_M\, (X,a,b) \neq 0$ and ${\rm corank}_M \,(J \setminus \{a,b\})
= 0$. It follows from Lemmas \ref{dual flats}, \ref{isthmus 2} that
$\hat{X}(a,b)$ is a flat in $M'$.

Clearly, $a,b \not \in \hat{X} (a,b)$ and ${\rm rank}_{M'} \{ a,b \}
= 2$.

Since $(X,a,b)$ is a parallelism, ${\rm rank}_M X \cup \{a,b \} =
{\rm rank}_M X +1$. So, ${\rm corank}_M X = {\rm corank}_M X \cup \{
a,b\} +1$. Thus, ${\rm nullity}_{M'} \hat{X} (a,b) \cup \{a,b \} =
{\rm nullity}_{M'} \hat{X} (a,b) + 1$. Hence ${\rm rank}_{M'}
\hat{X} (a,b) \cup \{ a,b\} = {\rm rank}_{M'} \hat{X} (a,b) +1$.
This proves that $(\hat{X} (a,b), a,b)$ is a parallelism in $M'$.

Furthermore, ${\rm l}_{M} X = b^{10}_{M/X} = b^{10}_{M' - X} = {\rm
w}_{M'} (\hat{X}(a,b),a,b)$ by Lemmas \ref{b10cor}, \ref{delcont}.
Hence, ${\rm vol}_M (X,a,b) = {\rm vol}_{M'} (\hat{X}(a,b),a,b)$.
\end{proof}

\subsection{Proof of Lemma \ref{product equals negative one}}
\label{proof of product equals negative one}
\begin{lemma} \label{vertex on parallel edge}
Let $(L , H^a , H^{b})$ be a parallelism in $A[\tau]$ and let
$P=H^{j_1}\cap \cdots \cap H^{j_k}$ be a vertex on $L \setminus
H^{b}$ for some $I=\{j_1, \ldots , j_k\} \subset J$. Let
$\hat{I}(a,b)= J \setminus (I \cup \{a,b \})$. Then
$|\hat{I}(a,b)|=n$ and $\check{P} = \cap_{j \in \hat{I}(a,b)} H_j$
is a vertex on $\hat L(H^a,H^{b}) \setminus H_{b}$ in the dual
arrangement $A[\tau']$.
\end{lemma}

\begin{proof}
Since ${\rm corank}_{M[\tau]} I =1$ and $a \not \in X$, ${\rm
rank}_{M[\tau]} I \cup \{ a\} = {\rm rank}\, M$. Thus, ${\rm
nullity}_{M[\tau]} I \cup \{a,b \} =1$ and hence ${\rm
corank}_{M[\tau]} \hat{I} (a,b) =1$. So, $\check P$ is a vertex on
$\hat{L} (H^a , H^b)$.

It remains to show that $\check{P} \not \in H_b$. Assume that
$\check{P} \in H_b$. Then ${\rm nullity}_{M[\tau']} \hat{I} (a,b) +
1 = {\rm nullity}_{M[\tau']} \hat{I} (a,b) \cup \{ b\}$. Thus, ${\rm
corank}_{M[\tau]} I \cup \{a,b\} + 1 = {\rm corank}_{M[\tau]} I \cup
\{ a\}$. But, ${\rm corank}_{M[\tau]} I \cup \{ a\} =0$. This is a
contradiction.
\end{proof}

\begin{proof}[Proof of Lemma \ref{product equals negative one}]
Let $(L , H^j , H^{N+1})$ be a parallelism in $A[\tau]$ and
\linebreak $(\hat{L}(H^j , H^{N+1}) , H_j , H_{N+1})$ the dual
parallelism in $A[\tau']$. Let $P$ and $\check{P}$ be vertices on $L
\setminus H^{N+1}$ and $\hat{L}(H^j , H^{N+1}) \setminus H_{N+1}$,
respectively, as in Lemma \ref{vertex on parallel edge}. We want to
show that $f^j (P) \cdot f_j (\check{P}) = -1$.

Let $I$ be as in Lemma \ref{vertex on parallel edge}. Assume that $I
= \{j_1 , j_2 , \ldots , j_k \}$, $j_1 < j_2 < \cdots < j_k$, and
$\hat{I}(j , N+1) = \{ j_{k+1}, \ldots , j_{k+n}\}$, $j_{k+1}<
\cdots < j_{k+n}$.  Let $m$ and $m'$ be the number of elements less
than $j$ in $I$ and $\hat{I}(a,b)$, respectively. Then $m + m' =
j-1$.

By Lemma \ref{value of fj},
$$
f^j (P) = (-1)^{k + m  }  \, \frac{B [j_1 , \dots , j_m, j ,
j_{m+1}, \ldots , j_k]}{B [j_1 , \dots , j_k , N+1]}
$$
and
$$
f_j (\check{P}) = (-1)^{n + m'  }  \, \frac{C [j_{k+1} , \dots ,
j_{k+m'}, j , j_{k+m'+1}, \ldots , j_{k+n}]}{C [j_{k+1} , \dots ,
j_{k+n} , N+1]}\ .
$$
By Theorem \ref{minors}, $f^j (P) \cdot f_j (\check{P}) =
(-1)^{(j+N+1) + (k+m)+ (n+m')} = -1 $.
\end{proof}

\subsection{Proof of Theorem \ref{main theorem}} \label{proof of main theorem}

\begin{lemma} \label{sum is beta}
Let $M$ and $M'$ be dual matroids. Let $X$ be a flat in $M$, $|X|=1$.
\begin{enumerate}
\item[$\bullet$] If $X$ is also a flat in $M'$, then ${\rm l}_M X + {\rm l}_{M'}X =
b^{10}_M$.
\item[$\bullet$] If $X$ is not a flat in $M'$, then ${\rm l}_M X =
b^{10}_M$.
\end{enumerate}
\end{lemma}
\begin{proof}
Suppose that $X$ is a flat in both $M$ and $M'$. Then ${\rm l}_M X =
b^{10}_{M/X}$ and ${\rm l}_{M'} X = b^{10}_{M'/X} = b^{10}_{M-X}$.
Hence, ${\rm l}_M X + {\rm l}_{M'}X =b^{10}_{M/X}+
b^{10}_{M-X}=b^{10}_M$, by Theorem \ref{tuttepolynomial}.

Let $X$ be a flat in $M$. Suppose that $X$ is not a flat in $M'$.
Then the matroid $M-X$ has an isthmus by Lemma \ref{dual flats}.
Hence, $b^{10}_{M-X} =0$ by Theorem \ref{tuttepolynomial}. Hence,
$b^{10}_{M/X} =b^{10}_M$, by Theorem \ref{tuttepolynomial}. That is,
${\rm l}_M X = b^{10}_M$.
\end{proof}

\begin{theorem} \label{betaprod}
Let $( \mathcal{A} [\tau] , \alpha )$ and $( \mathcal{A}
[\tau^{\prime}] , \alpha )$ be dual weighted affine arrangements.
Then
$$
{\rm B} ( \mathcal{A} [\tau] , \alpha ) \cdot {\rm B} ( \mathcal{A}
[\tau^{\prime}] , \alpha ) = \left[ \frac{\prod_{j \in \mathcal{J} }
\Gamma (\alpha_j + 1)}{\Gamma ( \sum_{ j \in \mathcal{J}} \alpha_j +
1 )} \right]^{\beta} \, .
$$
\end{theorem}
\begin{proof}
Let $L$ be a spacious edge of $A[\tau]$ that is not a hyperplane. By
Lemma \ref{dual edges}, the dual edge $\hat L$ is spacious in
$A[\tau']$. Furthermore, ${\rm vol}_{A[\tau]} L = {\rm
vol}_{A[\tau']} \hat{L}$.

Clearly, if $L \in \mathcal{L}_- [\tau]$, then $\hat{L} \in
\mathcal{L}_+ [\tau']$ and if $L \in \mathcal{L}_+ [\tau]$, then
$\hat{L} \in \mathcal{L}_- [\tau']$.
We also have $\alpha(L) + \alpha(\hat{L}) = 0$.

If    $L \in
\mathcal{L}_+ [\tau]$, then
\bea
\Gamma (\alpha(L) + 1)^{{\rm vol}_{A[\tau]}
L} = \Gamma (- \alpha(\hat{L}) + 1)^{{\rm vol}_{A[\tau']} \hat{L}}\ .
\eea
If $L \in \mathcal{L}_- [\tau]$, then
\bea
\Gamma (-
\alpha(L) + 1)^{{\rm vol}_{A[\tau]} L} = \Gamma ( \alpha(\hat{L}) +
1)^{{\rm vol}_{A[\tau']} \hat{L}}\ .
\eea

Hence, by Lemma \ref{sum is beta},
\begin{eqnarray}
\nonumber {\rm B} ( \mathcal{A}[\tau] , \alpha ) \cdot {\rm B} (
\mathcal{A}[\tau'] , \alpha ) & = & \frac{\prod_{j = 1}^N \Gamma
(\alpha_j +1 )^{{\rm l}_{A[\tau]} (H^j)} \cdot \Gamma (\alpha_j +1
)^{{\rm l}_{A[\tau']} (H_j)}}{\Gamma (-\alpha_{N+1} + 1)^{{\rm
l}_{A[\tau]} (H^{N+1})} \cdot
 \Gamma (-\alpha_{N+1} + 1)^{{\rm l}_{A[\tau']} (H_{N+1})}}   \\
\noalign{\medskip} \nonumber & = & \left[ \frac{\prod_{j = 1}^N
\Gamma (\alpha_j + 1)}{ \Gamma ( \sum_{j=1}^N \alpha_j
+1)}\right]^{\beta}
\end{eqnarray}
as desired.
\end{proof}

\begin{proof}[Proof of Theorem \ref{main theorem}]
Let $(\hat{L} (H^j , H^{N+1}) , H_j , H_{N+1})$ be the parallelism
in $A[\tau']$ dual to the parallelism $(L , H^j , H^{N+1})$ in $A
[\tau]$. Then
$$
 c \,({(f^j)}^{\alpha_j} , L )\, \cdot \,
 c \, ( (f_j)^{\alpha_j} , \hat{L} (H^j , H^{N+1})) \,
  = e^{i \pi \alpha_j} .
  $$
For each $j \in \mathcal{J}$, there are $\beta$ critical values of
$(f^j)^{\alpha_j}$. Hence
$$
\prod_{\stackrel{j \in \mathcal{J}}{\Delta \in \, {\tt Ch
\mathcal{A} [\tau]}}} c \,({(f^j)}^{\alpha_j} , \Delta ) \cdot
\prod_{\stackrel{j \in \mathcal{J}}{\Delta' \in \, {\tt Ch
\mathcal{A} [\tau']}}} c \,({(f_j)}^{\alpha_j} , \Delta' ) = e^{i
\pi \beta \sum_{j \in \mathcal{J}} \alpha_j} \ .
$$
The theorem now follows from Theorem \ref{betaprod}
\end{proof}

\section{Appendix A: Weak duality }
\label{Appendix A}

\subsection{Statement of results}
Let $\X$, $\W$ and $\X'$, $\W'$ denote the same spaces as in Section \ref{admissible
pairs}. Assume that for any $j \in J$, the functions $e^j \vert_{\W}$ and
$e_j \vert_{\W'}$ are not identically zero.
The pair $\tau=(\X,\W)$ with this property will be called a {\it
weakly admissible} pair in $\X$. Similarly, the pair $\tau'=(\X',
\W')$ with this property will be called a {\it weakly admissible}
pair in $\X'$. The pairs $\tau$ and $\tau'$ will be called {\it
weakly dual}.

Clearly, any dual pairs $\tau$ and $\tau'$ are weakly dual.

Let $\tau=(\X , \W)$ be a weakly admissible pair. As in Section
\ref{duality}, define the following objects. For $j \in J$, denote
$E^j = \{\, x \in \W\ :\ e^j (x) = 0\,\}$. These are vector subspaces of
$\W$ of codimension one.

Denote
\begin{enumerate}
\item [$\bullet$] $V = \Pee (\W)$, the projective space of dimension
$k$,
\item [$\bullet$]  $H^j = \Pee (E^j)$, $j \in J$, projective
hyperplanes in $V$,
\item [$\bullet$]
 $A [\tau] = \{\,H^j\ :\  j \in J \,\}$, the arrangement of
projective hyperplanes in $V$.
\end{enumerate}

For weakly dual pairs $\tau$ and $\tau'$, the corresponding projective arrangements
$A[\tau]$ and $A[\tau']$  will be called {\it weakly dual}.

For a weakly admissible pair $\tau$, introduce the matroid of $\tau$,
denoted $M[\tau]$, as the matroid of the collection of vectors
$\{\,e^j\vert_{\W}\ :\ j\in J \,\}$ in the dual space of $\W$.

Let $\tau$ and $\tau'$ be weakly dual pairs. Then the matroids
$M[\tau]$ and $M[\tau']$ are dual.

\begin{theorem} \label{conjecture1}
Let $A[\tau]$ and $A[\tau']$ be dual projective arrangements. Let L be a
spacious edge of $A[\tau]$ and $\hat L$ the dual spacious edge of
$A[\tau']$. Then the induced arrangement on $L$, denoted by
$A[\tau]_L$, is weakly dual to the
projective localization of $A[\tau']$ at $\hat L$, denoted by
$P(A[\tau']^{\hat L})$.
\end{theorem}

We will prove a more general result.
To formulate the result we need the notion of the
projective localization of a sub-arrangement.

\medskip

Let $A = \{\,H^j\ :\ j \in J \,\}$ be a projective arrangement of hyperplanes.
For a subset  $I \subset J$, the  set of hyperplanes $\{\,H^j\ :\ j\in I \,\}$
will be called {\it a sub-arrangement} of $A$.

Let $I \subset J$ be such that $\cap_{j \in I} H^j \neq \emptyset$. Then the
projectivization of the
central arrangement $\{\,H^j\ :\ j\in I \,\}$
will be called {\it the projective localization of the sub-arrangement}
$\{\,H^j\ :\ j\in I \,\}$, see Section \ref{arrangement and edges}.

\begin{theorem}\label{dual induced arrangements}
Let $A[\tau] = \{\,H^j\ :\ j \in J\, \}$ and $A[\tau'] =
\{\,H_j\ :\ j\in J\, \}$ be
weakly dual arrangements. Let $L$ be an edge of $A[\tau]$ and let
$X = \{\,j\in J\ :\ L \subset H^j \,\}$ be the flat associated to $L$ in the
matroid $M[\tau]$. Assume that $|X| < |J|$.

Let \,$\cap_{j\in \hat{X}}H_j = \emptyset$. Then the induced
arrangement on $L$, $A[\tau]_L$, is weakly dual to the sub-arrangement
$\{\,H_j\ :\ j\in \hat{X}\, \}$ of $A[\tau']$.

Let \,$\cap_{j\in \hat{X}}H_j \neq \emptyset$. Then the induced arrangement on $L$,
$A[\tau]_L$, is weakly dual to the projective localization of the
sub-arrangement $\{\,H_j\ :\ j\in \hat{X} \,\}$ of $A[\tau']$.
\end{theorem}

\subsection{Proof of Theorem \ref{dual induced arrangements}}
Let $\tau=(\X , \W)$ and $\tau'=(\X' , \W')$ be the weakly dual
admissible tuples.

For $j \in J$, denote $F^j = \{e^j=0 \} \subset \X$ and $F_j = \{
e_j=0\} \subset \X'$.

Consider the vector space $\cap_{j \in X} F^j$ with the basis $e_j$,
$j \in \hat{X}$. Then $\X' / \cap_{j \in \hat{X}} F_j$ may be
identified with the dual space of $\cap_{j \in X} F^j $ with the
dual basis $e^j + \cap_{j \in \hat{X}} F_j$, $j \in \hat{X}$.

The classes of the elements of $\W'$ in the vector space $\X' /
\cap_{j \in \hat{X}} F_j$ form a subspace of $\X' / \cap_{j \in
\hat{X}} F_j$. Denote this subspace by $\W' / \cap_{j \in \hat{X}}
F_j$.

Consider the tuples
$$
\sigma = \left( \cap_{j \in X} F^j \ , \ \cap_{j \in X} E^j \right)
, \quad   \sigma' = \left( \X' / \cap_{j \in \hat{X}} F_j \ ,  \ \W'
/ \cap_{j \in \hat{X}} F_j \right)  .
$$

\noindent 
{\bf Remark.} To make sense of the definitions of $\sigma$
and $\sigma'$ the following identifications are required.
\begin{enumerate}
\item[$\bullet$] The element $e^j + \cap_{j \in \hat{X}} F_j \in \X' /
\cap_{j \in \hat{X}} F_j $ defines a linear function on $\cap_{j \in
X} F^j $ given by $(e^j + \cap_{j \in \hat{X}} F_j) (x) = e^j (x)$ for
every $x \in \cap_{j \in X} F^j$.
\item[$\bullet$] The element $e_j \in \cap_{j \in \hat{X}} F_j$, $j \in
\hat{X}$, defines a linear function on $\X' / \cap_{j \in \hat{X}}
F_j$ given by $e_j (x' + \cap_{j \in \hat{X}} F_j) = e_j (x')$, for
every $x' \in \X'$.
\end{enumerate}

\begin{lemma}\label{sigmas}
The tuples $\sigma$ and $\sigma'$ are weakly dual.
\end{lemma}
\begin{proof}
We first show that $\sigma$ and $\sigma'$ are weakly admissible. For
$i \in \hat{X}$,
\linebreak
 $(e^i + \cap_{j \in \hat{X}} F_j) \vert_{\cap_{j
\in X} E^j} \equiv 0$  implies that $e^i \vert_{\cap_{j \in X} E^j }
\equiv 0$. Thus, $\cap_{j \in X} E^j \subset E^i$. This contradicts
the assumption that $X$ is a flat.

For any $i\in \hat{X}$, $e_i \vert_{\W' / \cap_{j \in \hat{X}} F_j}$
is a nonzero function since  $e_i \vert_{\W'}$ is a nonzero linear
function.  This shows that $\sigma$ and $\sigma'$ are weakly
admissible tuples.

The space  $\W' / \cap_{j \in \hat{X}} F_j$
clearly annihilates $\cap_{j \in X} E^j$.

We have $\dim \,\cap_{j \in X} E^j = \dim \,\W \cap(\cap_{j \in X}
F^j) =  {\rm corank}_{M[\tau]} X$. Hence the annihilator
of $\cap_{j \in X} E^j$  in $\X' /
\cap_{j \in \hat{X}} F_j$ has dimension $|\hat{X}| - {\rm
corank}_{M[\tau]} X = |\hat{X}| - {\rm nullity}_{M[\tau']} \hat{X} =
{\rm rank}_{M[\tau']} \hat{X}$.

Observe that the spaces $\W' / \cap_{j \in \hat{X}} F_j $ and $\W' /
\cap_{j \in \hat{X}} E_j $ are isomorphic. Hence, $\dim \,\W' /
\cap_{j \in \hat{X}} F_j = \dim \,\W' - \dim \,\cap_{j \in \hat{X}}
E_j = {\rm rank \,} M[\tau'] - {\rm corank}_{M[\tau']} \hat{X} =
{\rm rank}_{M[\tau']} \hat{X}$.

This shows that $\W' / \cap_{j \in \hat{X}} F_j$ is the annihilator
of $\cap_{j \in X} E^j$ in $\X' / \cap_{j \in \hat{X}} F_j$.
Hence, $\sigma$ and $\sigma'$ are weakly dual.
\end{proof}

\begin{lemma} \label{the isomorphism}
${}$

\begin{enumerate}
\item [(i)]
The arrangement $A[\sigma]$ is the induced arrangement $A[\tau]_L$.
\item [(ii)]
If $\cap_{j\in \hat{X}}H_j = \emptyset$, then the arrangement
$A[\sigma']$ is the sub-arrangement \linebreak
 $\{ H_j\ :\ j\in \hat{X} \}$ of
$A[\tau']$.
\item [(iii)]
If $\cap_{j\in \hat{X}} H_j \neq \emptyset$, then the arrangement
$A[\sigma']$ is the projective localization of the sub-arrangement
$\{ H_j\ :\ j\in \hat{X} \}$
 of $A[\tau']$.
\end{enumerate}
\end{lemma}
\begin{proof}
Statement (i) is clear.

Let $\U \subset \W'$ be a subspace such that $\U \oplus \cap_{j \in
\hat{X}}E_j = \W'$. Consider the isomorphism
$$
\U\   \to \ \W' /\cap_{j \in \hat{X}} F_j \ , \qquad
 u \ \mapsto\  u +  \cap_{j \in \hat{X}} F_j\ .
$$
For every $j\in \hat{X}$, the subspace $E_j \cap \U \subset \U$
corresponds to the subspace $ E_j / \cap_{j \in \hat{X}} F_j \subset
\W' / \cap_{j \in \hat{X}} F_j$. The arrangement of hyperplanes
$\Pee (\U \cap E_j) \subset \Pee (\U)$, $j\in \hat{X}$,\ in the
projectivization $\Pee (\U)$ is the arrangement described in part
(ii), if  \ $\cap_{j\in \hat{X}}H_j = \emptyset$,\ and is the
arrangement described in part (iii), if $\cap_{j\in \hat{X}} H_j
\neq \emptyset$.\ Statements (ii) and (iii)  are proved.
\end{proof}

Theorem \ref{dual induced arrangements} is a corollary of
 Lemmas \ref{sigmas} and \ref{the isomorphism}.

\section{Appendix B: Plucker coordinates of dual arrangements }
\label{Appendix B}
We formulate a statement which helps to
determine if two given arrangements are dual.

Let $k, n, N$   denote the same natural numbers as in Section 
\ref{admissible pairs}. 

Let $P(k+1,N+1)$ be the real projective space of dimension ${N+1
\choose k+1}-1$, whose projective coordinates $(\la_L)$ are labeled by
subsets $L = (l_1,\dots, l_{k+1})$ such that $1\leq l_1 <\dots <
l_{k+1} \leq N+1$. 

Similarly, let $P(k+1,N+1)$ be the real projective
space of (the same) dimension ${N+1 \choose n+1}-1$, whose projective coordinates
$(\mu_M)$ are labeled by subsets $M = (m_1,\dots, m_{n+1})$ such that
$1\leq m_1 <\dots < m_{n+1} \leq N+1$. 

Let 
$$
\delta\ :\ P(k+1,N+1)\ \to\ P(n+1,N+1)\ , 
\quad
 \la\ \mapsto\ \mu = \delta(\la)\ ,
$$
be the isomorphism, where for any 
$M = (m_1,\dots, m_{n+1})$, we set $\mu_M = (-1)^\sigma \la_{L}$,
where $L = ( l_{1},\dots, l_{k+1})$ is the subset complementary
to $M$  in $\{1, \dots,N+1\}$, and $\sigma = 1 + 2 + \dots + k + 1 + 
l_{1} + \dots +l_{k+1}$.

Let   $\X$ and $\X'$ denote the same spaces as in Section 
\ref{admissible pairs}. 
Let ${\rm Gr}$ denote the Grassmannian of all $k+1$-dimensional vector
subspaces of $\X$ and ${\Gr}'$ the Grassmannian of all
$n+1$-dimensional vector subspaces of $\X'$.
Let 
$$
\pi\ :\ \Gr\ \to\ P(k+1,N+1)\ , 
\quad
\W \ \mapsto\ (\la_L)= \pi(\W)\ ,
$$ 
with
$(\la_L) =  ( e^{l_1}\wedge \dots \wedge e^{l_{k+1}}\vert_{\W})$,
and
$$
\pi'\ :\ \Gr'\ \to\ P(n+1,N+1)\ ,
\quad 
\W'\ \mapsto\ (\mu_M)= \pi'(\W')\ ,
$$ 
with
$(\mu_M) =  ( e_{m_1}\wedge \dots \wedge e_{m_{n+1}}\vert_{\W'})$,
 be the Plucker imbeddings.

\begin{lemma}
\label{Plucker}
For $\W \in \Gr$ and $\W' \in \Gr'$,\ the subspace 
$\W'$ is the annihilator of the subspace $\W$ if and only if
$\pi' (\W') = \delta ( \pi (\W))$.
\end{lemma}

The lemma follows from Lemma \ref{value of fj}.

\end{document}